\documentclass[11pt]{article}%
\usepackage{amssymb}
\usepackage{amsfonts}
\usepackage{amsmath}
\usepackage{authblk}
\usepackage{url}
\usepackage{hyperref}
\usepackage{wrapfig}
\usepackage{xcolor}
\usepackage{soul}
\usepackage{graphicx}%
\usepackage{mwe}
\usepackage{subcaption}
\usepackage{lineno}
\providecommand{\U}[1]{\protect\rule{.1in}{.1in}}
\newtheorem{theorem}{Theorem}

\newtheorem{definition}[theorem]{Definition}

\begin{document}

\title{Number of paths in a graph}

\author[1]{\small Ivan Joki\'{c}\thanks{\emph{I.Jokic@tudelft.nl}}}
\author[1]{\small Piet Van Mieghem}

\affil[1]{\footnotesize Faculty of Electrical Engineering, Mathematics, and Computer Science,
Delft University of Technology, P.O. Box 5031, 2600 GA Delft, The Netherlands}

\maketitle

\begin{abstract}
The $k$-th power of the adjacency matrix of a simple undirected graph represents the number of walks with length $k$ between pairs of nodes. As a walk where no node repeats, a path is a walk where each node is only visited once. The set of paths constitutes a relatively small subset of all possible walks.
We introduce three types of walks, representing subsets of all possible walks. Considered types of walks allow for deriving an analytic solution for the number of paths of a certain length between node pairs in a matrix form. Depending on the path length, different approaches possess the lowest computational complexity.
We also propose a recursive algorithm for determining all paths in a graph, which can be generalised to directed (un)weighted networks.
\end{abstract}

\section{Introduction}\label{Sec_Intro}
Networks emerge naturally in many real-world systems \cite{barabasi2013network,newman2018networks}. From a Network Science perspective, a network consists of an underlying topology, called the graph, and a dynamic process taking place on the graph. 
Examples of real-world networks include infrastructural networks (such as road traffic networks, the Internet and the power grid networks), biological networks (such as the protein interaction network) and social networks.

A walk in a graph represents a sequence of nodes, where adjacent nodes in the sequence share a common link, as defined below in Definition \ref{Def_Walk}. A path, defined in Definition \ref{Def_Paths}, is a walk where no node repeats in the sequence.
Paths reflect the cost of connections between node pairs in a network, where the cost is quantified by weights of the links in the paths. If all link weights are constant or unity, then the weight of a path equals the hopcount \cite[Chapter 16]{van2014performance}, the number of hops or links in the path. A significant part of the literature on paths in graphs considers the problem of determining the number of paths, both in random graphs \cite{van2001paths} or in special types of graphs \cite{grzesik2022maximum}. However, to the best of our knowledge, an explicit solution for the number of paths of a certain length in a matrix form is not yet known. Our idea is to derive an analytic solution for the number of length $k$ paths between node pairs by removing those walks traversing a node multiple times from all possible walks with hopcount $k$, contained in the $k$-th power of the adjacency matrix. 

The problem of determining whether there is a path in a graph with at least $k$ links is NP-complete \cite{leiserson1994introduction}. Williams proposed in \cite{williams2009finding} an algorithm for finding paths with hopcount $k$,
while Schmid \textit{et al.} proposed an logarithm in \cite{schmid2017computing} for computing Tutte Paths. A graph is Hamiltonian if there is a path traversing each node in a graph. Bjorklund proposed in \cite{bjorklund2014determinant} a Monte Carlo algorithm that solves the Hamiltonian problem. Bax introduced an algorithm for the Hamiltonian path problem in \cite{bax1993inclusion}, based on the inclusion-exclusion formula. We generalise the approach of Bax in \cite{bax1993inclusion} and derive the solution for the number of paths with any hopcount $k$.

Section \ref{Sec_Notation} introduces the notion of walks and paths. We firstly introduce walks traversing a node multiple times in Section \ref{Sec_Repeating_Node} and derive an analytic solution for the number of paths between node pairs with hopcount $k\leq 4$. Section \ref{Sec_Walks_Node_Appearing_Once} analyses walks traversing a node exactly once, while those walks not traversing a node we examine in Section \ref{Sec_Walks_Node_Not_Appearing}. In Section \ref{Sec_Recursive_Algorithm}, we provide a recursive algorithm that identifies all possible paths in a graph. Finally, we conclude in Section \ref{Sec_Conclusion}.

\section{Walks in a graph}\label{Sec_Notation}
A graph $G(\mathcal{N},\mathcal{L})$ consists of $N = |\mathcal{N}|$ nodes, interconnected with $L=|\mathcal{L}|$ links. The graph topology is defined by the $N\times N$ adjacency matrix $A$, with matrix element $a_{ij} = 1$ if there is a link between node $i$ and node $j$, otherwise $a_{ij} = 0$. The degree vector $d=A\cdot u$, where $u$ is the $N\times 1$ all-one vector, contains the degree of each node and the degree $d_{i} = \sum_{j=1}^{N}a_{ij}$ equals the number of direct neighbors of node $i$. The $N\times N$ degree diagonal matrix $\Delta = \text{diag}(d)$ contains the degree of the nodes on its main diagonal.

\begin{definition}\label{Def_Walk}
	A walk of length $k$ from node $i$ to node $j$ is a sequence of $k$ links of the form $(n_{0}\rightarrow n_{1})(n_{1}\rightarrow n_{2})\dots (n_{k-1}\rightarrow n_{k})$, where $n_{0}=i$ and $n_{k} = j$.
\end{definition}
The length of a walk is the number of links in the walk and is often referred to as the walk hopcount \cite{van2010graph}. The first node in the sequence $n_{0}$ is the source node, while a walk ends with the destination node $n_{k}$. The Definition \ref{Def_Walk} naturally leads to the question of how many ways are there to reach node $i$ from node $j$, in $k$ hops.
\begin{theorem}\label{Theorem_Walks}
	The number of walks of length $k$ from node $i$ to node $j$ is equal to the element $\left(A^{k}\right)_{ij}$.
\end{theorem}
\textit{Proof} Provided in \cite[p. 26]{van2010graph} $\hfill \square $

The $k$-th power of the adjacency matrix $A^{k}$ contains the number of walks of length $k$ between each pairs of nodes in the graph. Any matrix derived solely from the adjacency matrix $A$, either via matrix product or the Hadamard product, carries information about walks. We denote the set of all possible walks of length $k$ as the $k$-dimensional walk space $\mathcal{W}[k]$, while the corresponding $N\times N$ walk matrix is $A^{k}$.
Since a matrix commutes with itself, it holds that
\begin{equation}\label{Eq_Spliting_Walks}
	A^{k} = A^{k-p}\cdot A^{p} =  A^{p}\cdot A^{k-p},
\end{equation}
where $p$ is here an integer between $1\leq p\leq k$, although (\ref{Eq_Spliting_Walks}) holds for any complex $p$ for which the matrix $A^p$ exists.  
The relation (\ref{Eq_Spliting_Walks}) teaches us that
a walk can be split into sub-walks, while the walk matrix can be obtained as the multiplication of walk matrices of the corresponding sub-walks.
\subsection{Paths in a graph}\label{Sec_Paths_Definition}
\begin{definition}\label{Def_Paths}
	A path is a walk in which all nodes are different. A path of length $k$ is defined by a sequence of $k+1$ node pairs: $(n_{0}\rightarrow n_{1})(n_{1}\rightarrow n_{2})\dots (n_{k-1}\rightarrow n_{k})$, where $n_{l}\neq n_{m}$ for all $0\leq l \neq m \leq k$.
\end{definition}
Paths account for a relatively small subset of all possible walks $\mathcal{W}[k]$ of length $k$. We denote the set of all possible paths of length $k$ as $\mathcal{P}[k]$, where $\mathcal{P}[k] \subseteq \mathcal{W}[k]$, while equality holds only for $k=1$, as shown later in (\ref{Eq_P_1}). The $N\times N$ path matrix $P_{k}$ contains the number of paths of length $k$ between any pair of nodes, with $(P_{k})_{ij}$ denoting the number of paths between node $i$ and node $j$, of hopcount $k$.

In the following sections, we introduce three different types of walk sets: walks with a node reappearing in the node sequence, walks where a node is not traversed and those walks traversing a node exactly once.
Based on each mentioned type of walks, we derive an analytic solution for the $N\times N$ path matrix $P_k$ with hopcount $k$.

\section{Node reappearance in a walk}\label{Sec_Repeating_Node}
A walk, introduced in Definition \ref{Def_Walk}, can traverse the same node multiple times, defining the first walk type we consider in this section.

\begin{definition}\label{Def_Repeating_Node_Walks}
	The set of all possible walks of length $k$, where the same node appears at least twice in the node sequence, on positions $i$ and $j$ (i.e. $n_{i} = n_{j}$), where $j>i+1$, is denoted as $\mathcal{W}_{(i,j)}[k]$. The $N\times N$ walk matrix with the number of such walks between any pair of nodes is denoted by $M(\mathcal{W}_{(i,j)}[k])$.
\end{definition}

\begin{figure}[!h]
	\begin{center}
		\includegraphics[ angle =0, scale= 0.58]{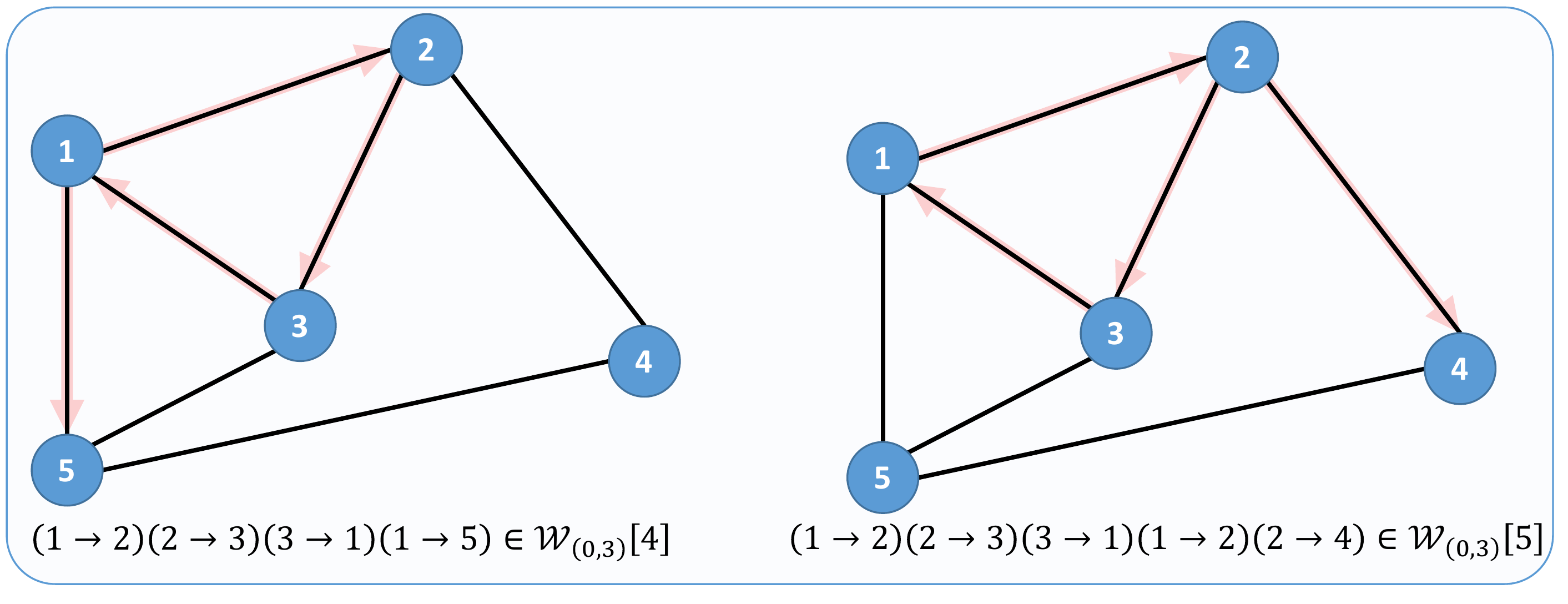}
		\caption{Examples of walks with node reappearance. Traversed links are colored in red, while arrows follow labeling in the node sequence.}
		\label{Fig_Node_reappearance}
	\end{center}
\end{figure}

From Definition \ref{Def_Repeating_Node_Walks}, we observe the following identity $\mathcal{W}_{(i,j)}[k] = \mathcal{W}_{(j,i)}[k]$, because $n_{i} = n_{j}$. The introduced constraint $j>i+1$ excludes two trivial cases. When $i=j$, the corresponding set of walks is actually the set of all possible walks of length $k$, i.e. $\mathcal{W}_{(i,i)}[k] = \mathcal{W}[k]$. On the other side, when $j=i+1$, the corresponding set of walks $\mathcal{W}_{(i,i+1)}[k] = \emptyset$, because a node cannot be adjacent to itself in a walk sequence, because there are no self-loops in a simple network. Figure \ref{Fig_Node_reappearance} illustrates a few examples of walks where a node is traversed multiple times.

It is of particular interest to consider walks, where the source node is also the destination node, i.e. the set $\mathcal{W}_{(0,k)}[k]$.

\begin{definition}\label{Def_Walk_Closed}
	A closed walk of length $k$ is a walk that starts in node $i$ and returns, after $k$ hops, to that same node $i$ (i.e. where $n_{0} = n_{k}$). The set of all possible closed walks of length $k$ is $\mathcal{W}_{(0,k)}[k]$. The corresponding $N\times N$ walk matrix is $M\left(\mathcal{W}_{(0,k)}[k]\right) = I\circ A^{k}$, where $I$ is the $N\times N$ identity matrix, while $\circ$ denotes the Hadamard product. The total number of closed walks of length $k$ in a graph is $\text{trace}(A^{k})$ and equals $\text{trace}(A^{k}) = \sum_{i=1}^{N}\lambda_{i}^{k}$, where $\lambda_j$ is the $j$-th largest eigenvalue of the adjacency matrix \cite{van2010graph}.
\end{definition}

\begin{figure}[!h]
	\begin{center}
		\includegraphics[ angle =0, scale= 0.62]{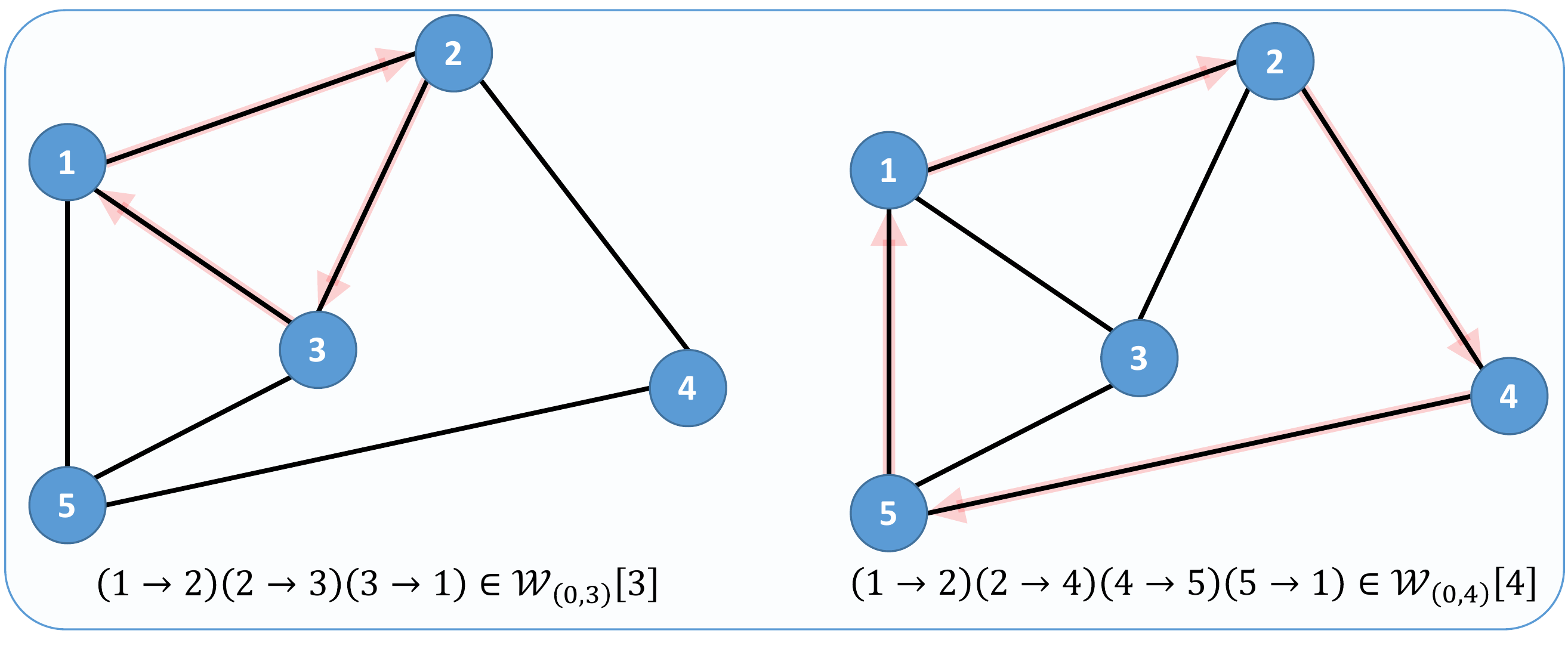}
		\caption{Examples of closed walks. Traversed links are colored in red, while arrows follow labeling in the node sequence.}
		\label{Fig_Closed_Walks}
	\end{center}
\end{figure}


A few examples of closed walks are depicted in Figure \ref{Fig_Closed_Walks}. A walk set $\mathcal{W}_{(i,j)}[k]$ can be split into three sub-walks: the set of all possible walks $\mathcal{W}[i]$ of length $i$, the set of closed walks $\mathcal{W}_{(0,j-i)}[j-i]$ of length $j-i$ and the set of all possible walks $\mathcal{W}[k-j]$ of length $k-j$. The corresponding $N\times N$ walk matrix $M\left(\mathcal{W}_{(i,j)}[k]\right)$, after applying Theorem \ref{Theorem_Walks}, Definition \ref{Def_Walk_Closed} and relation (\ref{Eq_Spliting_Walks}) becomes
\begin{equation}\label{Eq_Repeating_Node_Walk_Matrix}
	M\left(\mathcal{W}_{(i,j)}[k]\right) = A^{i}\cdot \left(I\circ A^{j-i}\right)\cdot A^{k-j}.
\end{equation}
The $N\times N$ walk matrix $M\left(\mathcal{W}_{(i,j)}[k]\right)$ in (\ref{Eq_Repeating_Node_Walk_Matrix}) is asymmetric, since the order of nodes in a walk sequence is labelled. In other words, if we label nodes of the walk sequence in $\mathcal{W}_{(i,j)}[k]$ from the destination node to the source node, we obtain the walk matrix $M\left(\mathcal{W}_{(k-j,k-i)}[k]\right) = A^{k-j}\cdot \left(I\circ A^{j-i}\right)\cdot A^{i}$, which together with (\ref{Eq_Repeating_Node_Walk_Matrix}) leads to the following identity
\begin{equation}\label{Eq_Repeating_Node_Transpose}
	M\left(\mathcal{W}_{(k-j,k-i)}[k]\right) = M^{T}\left(\mathcal{W}_{(i,j)}[k]\right).
\end{equation}
Identity (\ref{Eq_Repeating_Node_Transpose}) holds only for undirected networks, informing us that a walk can be performed from the source towards the destination node, but also in the reverse order, because all links are undirected.

\subsection{Analytic solution for the $N\times N$ path matrix $P_k$}\label{Sec_Paths}

\begin{theorem}\label{Theorem_All_Walks_and_paths}
	The set of all possible walks $\mathcal{W}[k]$ of length $k$ consists of the following subsets
	\begin{equation}\label{Eq_All_Walks_Theorem}
		\mathcal{W}[k] = \left(\bigcup_{i=0}^{k-2}\bigcup_{j=i+2}^{k}\mathcal{W}_{(i,j)}[k]\right) \cup \mathcal{P}[k].
	\end{equation}
\end{theorem}
\textit{Proof} A walk of length $k$ has either a repeating node in its sequence or represents a path of length $k$. Thus, by computing the set union of walk sets, defining walks with all possible repetitions of a node, and the set of paths, we obtain the set $\mathcal{W}[k]$ of all possible walks of length $k$ in (\ref{Eq_All_Walks_Theorem}), which completes the proof.  $\hfill \square $

A walk of length $k$ is either a path or a node is traversed multiple times. From Definition \ref{Def_Repeating_Node_Walks} and Definition \ref{Def_Paths}, we observe
\begin{equation}\label{Eq_Walks_Paths_Intersection}
	\mathcal{W}_{(i,j)}[k] \cap \mathcal{P}[k] = \emptyset,
\end{equation}
where $0\leq i \leq k-2$ and $j\geq i+2$. Based on Theorem \ref{Theorem_Walks}, equations (\ref{Eq_All_Walks_Theorem}) and (\ref{Eq_Walks_Paths_Intersection}), we obtain a general solution for the $N\times N$ path matrix $P_{k}$
\begin{equation}\label{Eq_Path_Matrix_General_Sets}
	P_{k} = A^{k} - M\left(\bigcup_{i=0}^{k-2}\bigcup_{j=i+2}^{k}\mathcal{W}_{(i,j)}[k]\right).
\end{equation}
The double set union in (\ref{Eq_Path_Matrix_General_Sets}) defines the union of $\frac{k\cdot(k-1)}{2}$ walk sets. Thus, the number of walk sets of the form $\mathcal{W}_{(i,j)}[k]$ increases as a square function of the hopcount $k$.
In general, the sets of walks $\mathcal{W}_{(i_{1},j_{1})}[k]$ and $\mathcal{W}_{(i_{2},j_{2})}[k]$ overlap, which complicates the computation of equation (\ref{Eq_Path_Matrix_General_Sets}). Therefore, we apply the inclusion-exclusion formula.

\subsection{Inclusion-exclusion  formula}\label{Sec_Inclusion_Exclusion}
The inclusion-exclusion formula \cite[p. 10]{van2014performance} defines the cardinality of the union of sets and thus transforms the second term on the right-hand side of the equation in (\ref{Eq_Path_Matrix_General_Sets}) as follows

\begin{equation}\label{Eq_Inclusion_Exclusion_Applied}
	\begin{split}
		M\left(\bigcup_{i=0}^{k-2}\bigcup_{j=i+2}^{k}\mathcal{W}_{(i,j)}[k]\right) = & 	\sum\limits_{i_{1}=0}^{k-2}\sum\limits_{j_{1}=i_{1}+2}^{k} M\left(\mathcal{W}_{(i_{1},j_{1})}[k]\right) \\
		- & \sum\limits_{i_{1}=0}^{k-2}\sum\limits_{j_{1}=i_{1}+2}^{k}\sum\limits_{i_{2}=i_{1}}^{k-2}\sum\limits_{j_{2}=q_2}^{k} M\left(\mathcal{W}_{(i_{1},j_{1})}[k] \cap \mathcal{W}_{(i_{2},j_{2})}[k]\right) \\
		+ & \dots + \\
		+ & (-1)^{k-1}\sum\limits_{i_{1}=0}^{k-2}\sum\limits_{j_{1}=i_{1}+2}^{k} \dots \sum\limits_{i_{k}=i_{k-1}}^{k-2}\sum\limits_{j_{k}=q_k}^{k} M\left(\bigcap_{z=1}^{k} \mathcal{W}_{(i_{z},j_{z})}[k]\right),
	\end{split}
\end{equation}
where $q_m = j_{m-1} + 2$ if $i_{m} = i_{m-1}$, otherwise $q_m = i_{m} + 2$ with $1<m\leq k$. Since there are $\frac{k\cdot(k-1)}{2}$ different walk sets with a repeating node, the total number of terms in the inclusion-exclusion formula in (\ref{Eq_Inclusion_Exclusion_Applied}) is $2^{\left(\frac{k\cdot(k-1)}{2}\right)}-1$. However, not all walk set intersections in (\ref{Eq_Inclusion_Exclusion_Applied}) define possible walks, because a node cannot be adjacent to itself in the walk sequence. Under the assumption that each matrix in (\ref{Eq_Path_Matrix_General_Sets}) has an analytic solution, the complexity of computing the $N\times N$ path matrix $P_k$ with the number of path of length $k$ between node pairs is $O\left(kN^3 2^{\frac{k^2-k}{2}} \right)$. 

\subsection{Recursive Solution for the $N\times N$ Path Matrix $P_{k}$}\label{Sec_Recursive_Solution}

In this subsection we reason why an analytic recursive solution for the $N\times N$ path matrix $P_{k}$ of the hopcount $k$ seems infeasible. From (\ref{Eq_Path_Matrix_General_Sets}), we derive the following identity
\[
P_{k} + M\left(\bigcup_{i=0}^{k-2}\bigcup_{j=i+2}^{k}\mathcal{W}_{(i,j)}[k]\right) = P_{k-1}\cdot A + M\left(\bigcup_{i=0}^{k-3}\bigcup_{j=i+2}^{k-1}\mathcal{W}_{(i,j)}[k-1]\right)
\]
from where we conclude that the $N\times N$ path matrix $P_{k}$ of length $k$ obeys the following recursion
\begin{equation}\label{Eq_Recursive_Solution_Main}
	P_{k} = P_{k-1}\cdot A - F_{k},
\end{equation}
where the $N\times N$ matrix $F_{k}$ is defined as
\begin{equation}\label{Eq_Recursive_Solution_Matrix_F}
	F_{k} = M\left(\bigcup_{i=0}^{k-2}\bigcup_{j=i+2}^{k}\mathcal{W}_{(i,j)}[k]\right) - M\left(\bigcup_{i=0}^{k-3}\bigcup_{j=i+2}^{k-1}\mathcal{W}_{(i,j)}[k-1]\right)\cdot A.
\end{equation}

In Appendix \ref{App_Theorem_Recursive}, we derive the first two sum terms of the $N\times N$ matrix $F_{k}$, that illustrate the difficulty to derive a complete closed form solution. To provide an argument for why the recursive solution does not seem possible, we denote a path $p_k=(n_0\rightarrow n_1)(n_1\rightarrow n_2)\dots (n_{k-1}\rightarrow n_k)$ as a node sequence, where $n_{l}\neq n_{m}$ for all $0\leq l \neq m \leq k$. The number of paths between node $i$ and node $j$ of length $k+1$ can be computed as follows
\begin{equation}\label{Eq_Paths_Condiitions}
    \left(P_{k+1}\right)_{ij} = \sum\limits_{p_k \in \mathcal{P}[k]} \mathbf{1}_{\{n_0=i\}} \mathbf{1}_{\{n_1 \neq j\}}\mathbf{1}_{\{n_2\neq j\}}\mathbf{1}_{\{n_3\neq j\}} \dots \mathbf{1}_{\{n_{k-1}\neq j\}},
\end{equation}
where the set of node $j$ neighbours is denoted as $\mathcal{N}_j$, while $\mathbf{1}_{x}$ is the indicator function that equals 1 if statement $x$ is true, otherwise $\mathbf{1}_{x}=0$.
On the other side, we observe from (\ref{Eq_Recursive_Solution_Main}) that the first term of the recursive solution
\[\left(P_{k+1}\right)_{ij} = \sum\limits_{m=1}^{N}   \left(P_{k}\right)_{im}\cdot a_{mj} -   \left(F_{k}\right)_{ij}\]
examines only the first two conditions from (\ref{Eq_Paths_Condiitions}), because $\sum_{p_k \in \mathcal{P}[k]}\mathbf{1}_{\{n_0 =i\}}\mathbf{1}_{\{n_k\in \mathcal{N}_j\}} = \left(P_k \cdot A\right)_{ij}$. The element $\left(P_k \cdot A \right)_{ij}$ contains the number of walks between node $i$ and node $j$ of length $k+1$, composed of all paths of length $k$, from node $i$ to an adjacent node $m\in \mathcal{N}_j$. However, not each such a walk represents a path. To obtain length $k+1$ paths, we need to distract from $\left(P_k \cdot A \right)_{ij}$ those paths of length $k$ between node $i$ and node $m \in \mathcal{N}_j$, traversing also node $j$. 
The $N\times N$ path matrix $P_k$ counts the number of length $k$ paths between pairs of nodes, without comprising information about traversed nodes per each path. 
Therefore, the recursive solution in (\ref{Eq_Recursive_Solution_Main}) requires manipulating an exponentially large number of walk matrices, to account for walks in $P_k\cdot A$ that are not paths.


In the following subsections, we derive an explicit form of the $N\times N$ path matrix $P_{k}$ up to hopcount $k \leq 4$. 

\subsection{Path matrix $P_1$ of length $k=1$}\label{Sec_P_1}

The $N\times N$ adjacency matrix $A$, by its definition, defines all paths of length $k=1$ and thus the $N\times N$ path matrix $P_{1}$
\begin{equation}\label{Eq_P_1}
	P_{1} = A,
\end{equation}
because there is a path of length $k=1$ between node $i$ and $j$ only if they share a link (i.e. if $a_{ij}=1$).
Only in case $k=1$, the set of all walks
\[
\mathcal{W}[1] =  \mathcal{P}[1]
\]
consists of only paths, because there are no self-loops in simple networks.

\subsection{Path matrix $P_2$ of length $k=2$}\label{Sec_P_2}
\begin{figure}[!h]
	\begin{center}
		\includegraphics[ angle =0, scale= 1.0]{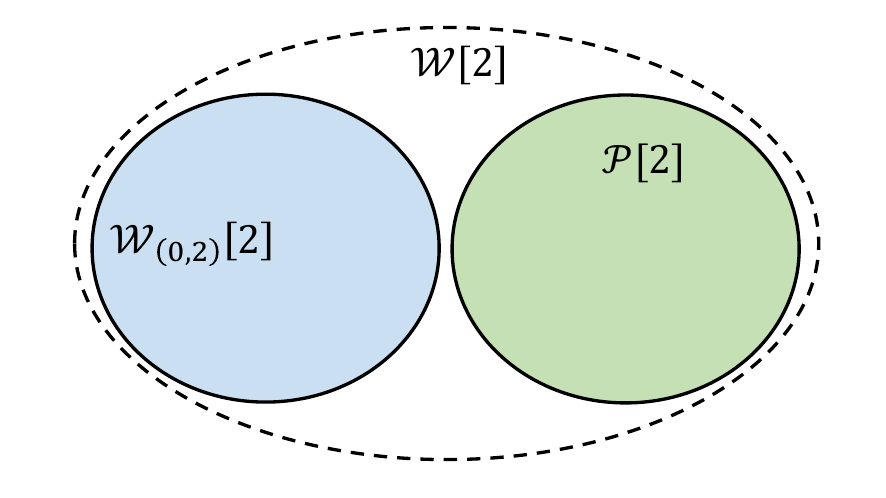}
		\caption{The set $\mathcal{W}[2]$ of all possible walks of length $k=2$ and subsets.}
		\label{Fig_Walks_Space_2}
	\end{center}
\end{figure}

When the hopcount $k>1$, walks with repeating nodes emerge, revealing the graph connectivity patterns. 
The only possible repetition of nodes in a walk with the hopcount $k=2$ is when $n_{0} = n_{2}$, forming the set $\mathcal{W}_{(0,2)}[2]$ of closed walks. Since walks of length $k=2$ consist of either closed walks or paths, the set $\mathcal{W}[2]$ of all possible walks of length $k=2$ is as follows
\[
\mathcal{W}[2] = \mathcal{W}_{(0,2)}[2] \cup \mathcal{P}[2].
\]
The Venn diagram for the walk space of the hopcount $k=2$ is provided on Figure \ref{Fig_Walks_Space_2}. From (\ref{Eq_Repeating_Node_Walk_Matrix}) we obtain the $N\times N$ walk matrix $M\left(\mathcal{W}_{(0,2)}[2]\right)$ as follows
\[
M\left(\mathcal{W}_{(0,2)}[2]\right) = I\circ A^{2}.
\]
By importing the above equation into (\ref{Eq_Path_Matrix_General_Sets}), we obtain the $N\times N$ path matrix $P_{2}$
\begin{equation}\label{Eq_P_2}
	P_{2} = A^{2} - I\circ A^{2}.
\end{equation}

\subsection{Path matrix $P_3$ of length $k=3$}\label{Sec_P_3}
The set $\mathcal{W}[3]$ of all walks with the hopcount $k=3$ represents the set union of the walk sets with any possible node repetition and the path set $\mathcal{P}[3]$. By applying (\ref{Eq_All_Walks_Theorem}) for the hopcount $k=3$, we obtain
\[
\mathcal{W}[3] = \mathcal{W}_{(0,2)}[3] \cup \mathcal{W}_{(0,3)}[3] \cup \mathcal{W}_{(1,3)}[3] \cup \mathcal{P}[3].
\]
By importing (\ref{Eq_Path_Matrix_General_Sets}), the above equation transforms
\begin{equation}\label{Eq_Walks_Number_k_3}
	\begin{split}
	A^{3} &= M\left(\mathcal{W}_{(0,2)}[3]\right) + M\left(\mathcal{W}_{(0,3)}[3]\right) + M\left(\mathcal{W}_{(1,3)}[3]\right) + P_{3} \\
	& - M\left(\mathcal{W}_{(0,2)}[3] \cap \mathcal{W}_{(0,3)}[3]\right) - M\left(\mathcal{W}_{(0,2)}[3] \cap \mathcal{W}_{(1,3)}[3]\right) - M\left(\mathcal{W}_{(0,3)}[3] \cap \mathcal{W}_{(1,3)}[3]\right) \\
	& + M\left(\mathcal{W}_{(0,2)}[3] \cap \mathcal{W}_{(0,3)}[3] \cap \mathcal{W}_{(1,3)}[3]\right). 
	\end{split}	  
\end{equation}

\begin{figure}[!h]
	\begin{center}
		\includegraphics[ angle =0, scale= 1.0]{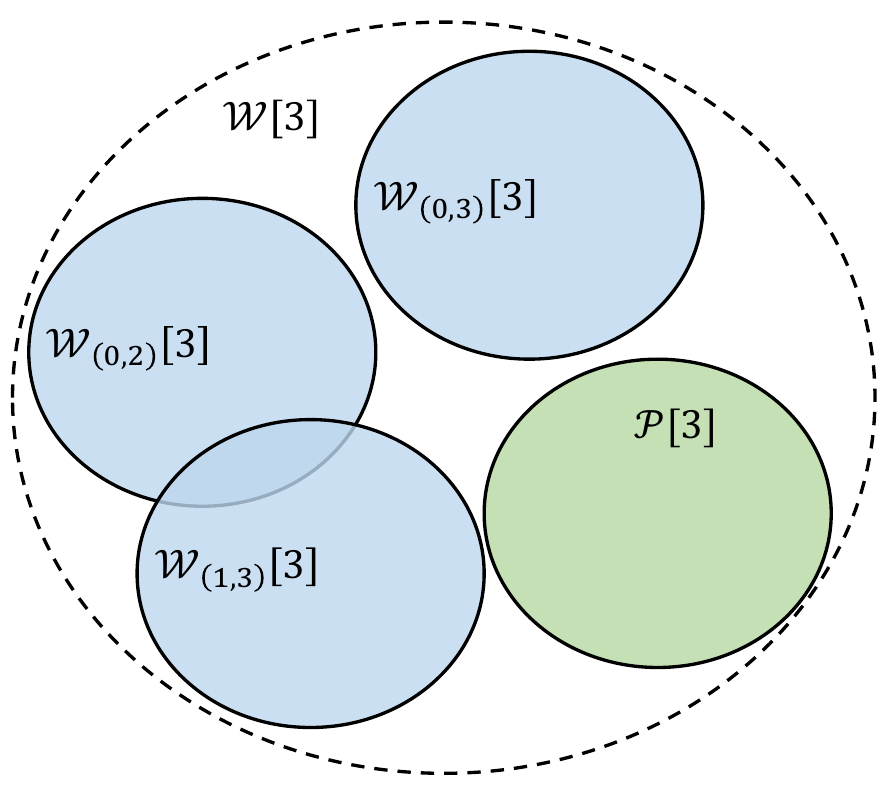}
		\caption{The set of all possible walks of length $k=3$, the subset of paths and subsets of walks with a repeating node.}
		\label{Fig_Walks_Space_3}
	\end{center}
\end{figure}

The set intersection $(\mathcal{W}_{(0,2)}[3]\cap\mathcal{W}_{(0,3)}[3])$ defines walks where $n_{0}=n_{2}=n_{3}$. Since a node is not adjacent to itself in a simple network (i.e. $a_{ii}=0$), such walks do not exist and thus $(\mathcal{W}_{(0,2)}[3]\cap\mathcal{W}_{(0,3)}[3]) = \emptyset$. The same reasoning holds for the sets $(\mathcal{W}_{(0,3)}[3]\cap\mathcal{W}_{(1,3)}[3]) = \emptyset$ and $(\mathcal{W}_{(0,2)}[3]\cap\mathcal{W}_{(0,3)}[3]\cap\mathcal{W}_{(1,3)}[3]) = \emptyset$, which simplifies the relation (\ref{Eq_Walks_Number_k_3})
\begin{equation}\label{Eq_Walks_Number_k_3_final}
		A^{3} =	M\left(\mathcal{W}_{(0,2)}[3]\right) + M\left(\mathcal{W}_{(0,3)}[3]\right) + M\left(\mathcal{W}_{(1,3)}[3]\right) + P_{3} - M\left(\mathcal{W}_{(0,2)}[3] \cap \mathcal{W}_{(1,3)}[3]\right). 
\end{equation}
The set $\mathcal{W}[3]$ of all possible walks with the hopcount $k=3$, the walk subsets $\mathcal{W}_{(0,2)}[3],\mathcal{W}_{(0,3)}[3],\mathcal{W}_{(1,3)}[3]$ with a repeating node and the path set $\mathcal{P}[3]$ are presented in Figure \ref{Fig_Walks_Space_3}. 

A walk of length $k=3$ where $n_{0}=n_{2}$ and $n_{1}=n_{3}$ starts from a node $i$, visits a neighbouring node $j$, traverses again the node $i$ and ends in the adjacent node $j$. Thus, for a pair of adjacent nodes $i$ and $j$, there is only one such a path. We denote the $N\times N$ corresponding walk matrix $M\left(\mathcal{W}_{(0,2)}[3] \cap \mathcal{W}_{(1,3)}[3]\right)$
\[
M\left(\mathcal{W}_{(0,2)}[3] \cap \mathcal{W}_{(1,3)}[3]\right) = A\circ A\circ A = A.
\]
Finally, after importing the above equation and (\ref{Eq_Repeating_Node_Walk_Matrix}) into (\ref{Eq_Walks_Number_k_3_final}), we obtain
\begin{equation}\label{Eq_Path_Matrix_3}
	P_{3} = A^{3} - \left(I\circ A^{2}\right) \cdot A - I\circ A^{3} - A\cdot \left(I\circ A^{2}\right) + A.
\end{equation}

\subsection{Path matrix $P_4$ of length $k=4$}\label{Sec_P_4}

\begin{figure}[!h]
	\begin{center}
		\includegraphics[ angle =0, scale= 0.9]{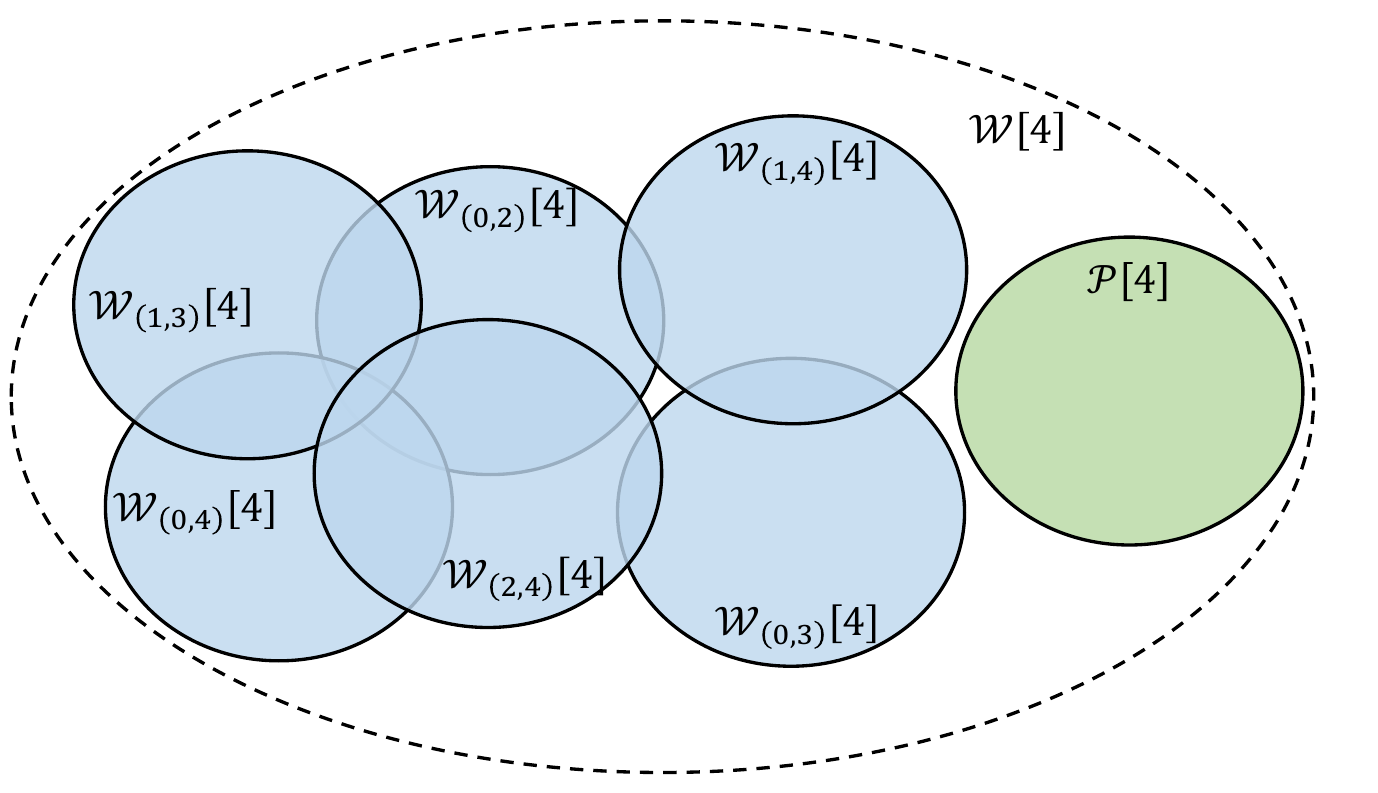}
		\caption{The set $\mathcal{W}[4]$ of all possible walks of length $k=4$, the subset of paths and subsets of walks with repeating nodes.}
		\label{Fig_Walks_Space_4}
	\end{center}
\end{figure}

The set $\mathcal{W}[4]$ of all possible walks of length $k=4$ represents the set union of the following sets
\[
\mathcal{W}[4] = \mathcal{W}_{(0,2)}[4] \cup \mathcal{W}_{(0,3)}[4] \cup \mathcal{W}_{(0,4)}[4] \cup \mathcal{W}_{(1,3)}[4] \cup \mathcal{W}_{(1,4)}[4] \cup \mathcal{W}_{(2,4)}[4] \cup P[4],
\]
The set $\mathcal{W}[4]$ of all walks with the hopcount $k=4$, the path set $\mathcal{P}[4]$ and the walk sets with a repeating node are presented in Figure \ref{Fig_Walks_Space_4}. 
Not all defined walk subsets with repeating nodes overlap, as presented in Figure \ref{Fig_Walks_Space_4}.
\begin{equation}\label{Eq_Walks_Number_k_4}
	\begin{split}
		A^{4} &= P_{3} + M\left(\mathcal{W}_{(0,2)}[4]\right) + M\left(\mathcal{W}_{(0,3)}[4]\right) + M\left(\mathcal{W}_{(0,4)}[4]\right) + M\left(\mathcal{W}_{(1,3)}[4]\right) + M\left(\mathcal{W}_{(1,4)}[4]\right) + M\left(\mathcal{W}_{(2,4)}[4]\right) \\
		& - M\left(\mathcal{W}_{(0,2)}[4] \cap \mathcal{W}_{(1,3)}[4]\right) - M\left(\mathcal{W}_{(0,2)}[4] \cap \mathcal{W}_{(1,4)}[4]\right) - M\left(\mathcal{W}_{(0,2)}[4] \cap \mathcal{W}_{(2,4)}[4]\right) \\
		& - M\left(\mathcal{W}_{(0,3)}[4] \cap \mathcal{W}_{(1,4)}[4]\right) - M\left(\mathcal{W}_{(0,3)}[4] \cap \mathcal{W}_{(2,4)}[4]\right) - M\left(\mathcal{W}_{(0,4)}[4] \cap \mathcal{W}_{(1,3)}[4]\right) \\
		&- M\left(\mathcal{W}_{(0,4)}[4] \cap \mathcal{W}_{(2,4)}[4]\right) - M\left(\mathcal{W}_{(1,3)}[4] \cap \mathcal{W}_{(2,4)}[4]\right) \\
		& + M\left(\mathcal{W}_{(0,2)}[4] \cap \mathcal{W}_{(1,3)}[4] \cap \mathcal{W}_{(2,4)}[4]\right) + M\left(\mathcal{W}_{(0,2)}[4] \cap \mathcal{W}_{(1,3)}[4] \cap \mathcal{W}_{(0,4)}[4]\right) \\
		& + M\left(\mathcal{W}_{(0,2)}[4] \cap \mathcal{W}_{(0,4)}[4] \cap \mathcal{W}_{(2,4)}[4]\right) + M\left(\mathcal{W}_{(0,4)}[4] \cap \mathcal{W}_{(1,3)}[4] \cap \mathcal{W}_{(2,4)}[4]\right) \\
		& - M\left(\mathcal{W}_{(0,2)}[4] \cap \mathcal{W}_{(0,4)}[4] \cap \mathcal{W}_{(1,3)}[4] \cap \mathcal{W}_{(2,4)}[4]\right). 
	\end{split}	  
\end{equation}
In the following part, we derive in sequel the walk matrices of the set intersections with three and four walk sets from the above relation.

Walk set $\left(\mathcal{W}_{(0,2)}[3] \cap \mathcal{W}_{(1,3)}[3] \cap \mathcal{W}_{(2,4)}[3]\right)$ defines walks with the hopcount $k=4$ where $n_{0} = n_{2} = n_{4}$ and $n_{1} = n_{3}$ originate from node $n_{0}$, visits node $n_{1}$, returns to node $n_{0}$ and repeats the same walk pattern, finishing at node $n_{0}$. 
\begin{equation}\label{Eq_M_k_4_02_13_24}
M\left(\mathcal{W}_{(0,2)}[4] \cap \mathcal{W}_{(1,3)}[4] \cap \mathcal{W}_{(2,4)}[4]\right) = I\circ \left(A\cdot \left(A\circ A\circ A\right) \right) = I\circ A^{2}.
\end{equation}

In addition, walk set $\left(\mathcal{W}_{(0,2)}[4] \cap \mathcal{W}_{(1,3)}[4] \cap \mathcal{W}_{(0,4)}[4]\right)$ and $\left(\mathcal{W}_{(0,4)}[4] \cap \mathcal{W}_{(1,3)}[4] \cap \mathcal{W}_{(2,4)}[4]\right)$ define walks of length $k=4$ where $n_{0}=n_{2}=n_{4}$ and $n_{1}=n_{3}$. Thus, we observe
\[
\left(\mathcal{W}_{(0,2)}[4] \cap \mathcal{W}_{(1,3)}[4] \cap \mathcal{W}_{(0,4)}[4]\right) = \left(\mathcal{W}_{(0,2)}[4] \cap \mathcal{W}_{(1,3)}[4] \cap \mathcal{W}_{(2,4)}[4]\right) = \left(\mathcal{W}_{(0,4)}[4] \cap \mathcal{W}_{(1,3)}[4] \cap \mathcal{W}_{(2,4)}[4]\right).
\]

On the contrary, the walk set $\left(\mathcal{W}_{(0,2)}[4] \cap \mathcal{W}_{(0,4)}[4] \cap \mathcal{W}_{(2,4)}[4]\right)$ defines walks with hopcount $k=4$ where $n_{0} = n_{2} = n_{4}$. Such walks start from a node $n_{0}$, visits an adjacent node, returns to node $n_{0}$, traverses an adjacent node once more and returns again to node $n_{0}$. The corresponding $N\times N$ walk matrix $M\left(\mathcal{W}_{(0,4)}[4] \cap \mathcal{W}_{(1,3)}[4] \cap \mathcal{W}_{(2,4)}[4]\right)$ is as follows
\begin{equation}\label{Eq_M_k_4_02_04_24}
	M\left(\mathcal{W}_{(0,2)}[4] \cap \mathcal{W}_{(0,4)}[4] \cap \mathcal{W}_{(2,4)}[4]\right) = \left(I\circ (A\cdot A)\right)\cdot \left(I\circ (A\cdot A)\right) = \left(I\circ A^{2}\right)^{2}.
\end{equation}

The set $ M\left(\mathcal{W}_{(0,2)}[4] \cap \mathcal{W}_{(0,4)}[4] \cap \mathcal{W}_{(1,3)}[4] \cap \mathcal{W}_{(2,4)}[4]\right)$ defines walks of length $k=4$ where $n_{0}=n_{2}=n_{4}$ and $n_{1}=n_{3}$ and thus
\[
M\left(\mathcal{W}_{(0,2)}[4] \cap \mathcal{W}_{(0,4)}[4] \cap \mathcal{W}_{(1,3)}[4] \cap \mathcal{W}_{(2,4)}[4]\right) = M\left(\mathcal{W}_{(0,2)}[4] \cap \mathcal{W}_{(1,3)}[4] \cap \mathcal{W}_{(2,4)}[4]\right) = I\circ A^{2}.
\]
Finally, we derive the $N\times N$ path matrix $P_{4}$ of length $k=4$
\begin{equation}\label{Eq_P_4}
	\begin{split}
		P_{4} =& A^{4} - \left(I\circ A^{2}\right)\cdot A^{2} - \left(I\circ A^{3}\right)\cdot A - I\circ A^{4} - A\cdot \left(I\circ A^{2}\right)\cdot A - A\cdot \left(I\circ A^{3}\right) - A^{2}\cdot \left(I\circ A^{2}\right) \\
		+ & 3\cdot \left(I\circ A^{2}\right)^{2} + 3\cdot A\circ A^{2} + I\circ \left(A\cdot \left(I\circ A^{2}\right) \cdot A\right) + 2\cdot A^{2} \\
		- & 3\cdot \left(I\circ A^{2}\right) - \left(I\circ A^{2}\right)^{2} \\
		+ & \left(I\circ A^{2}\right).
	\end{split}
\end{equation}

Determining walk matrices of all walk subsets in (\ref{Eq_Inclusion_Exclusion_Applied}) with increasing hopcount $k$ becomes intractable. 
While we provide above an explicit solution for the $N\times N$ path matrix $P_{k}$, with $k\leq 4$, already for $k=5$, providing the explicit enumeration of the number of paths between any pair of nodes becomes far more involving. 

\section{Walks traversing a node exactly once}\label{Sec_Walks_Node_Appearing_Once}

\begin{figure}[!h]
	\begin{center}
		\includegraphics[ angle =0, scale= 0.62]{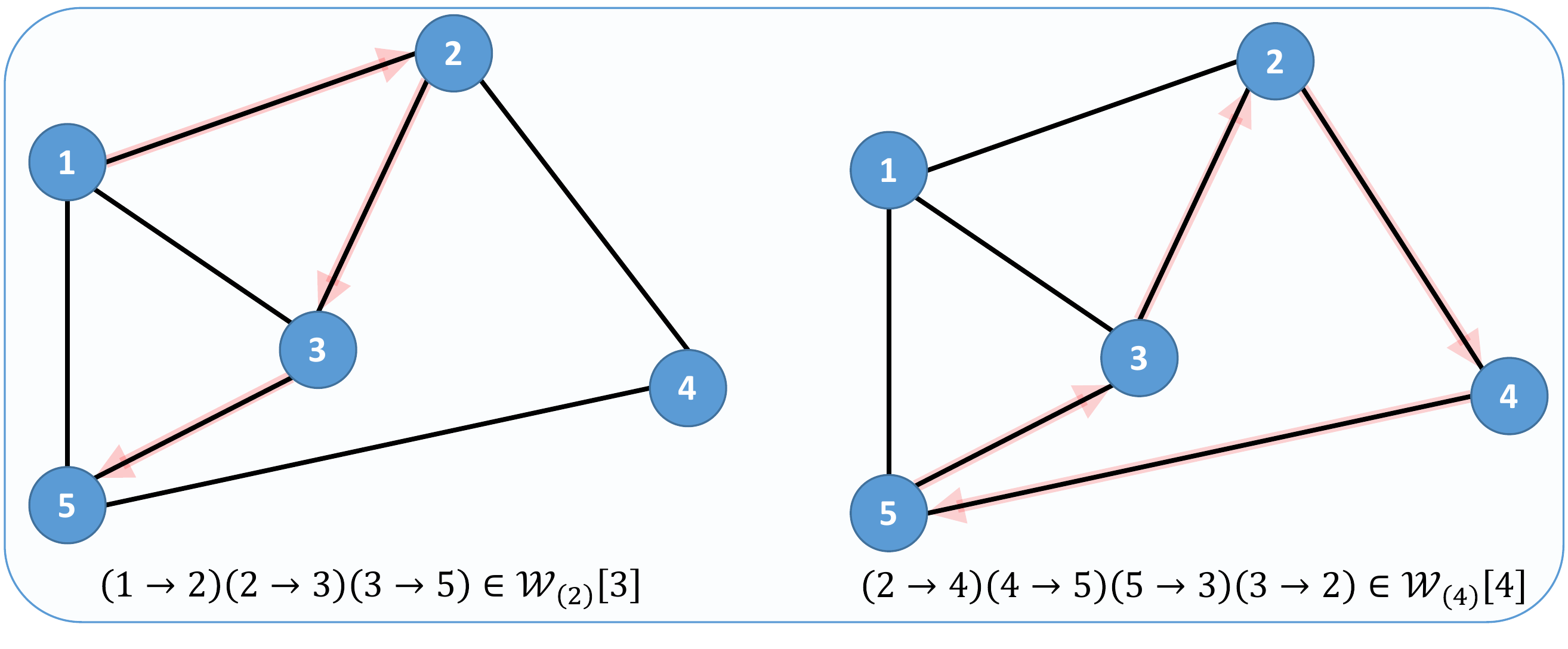}
		\caption{Examples of walks in which a node is traversed exactly once. Traversed links are colored in red, while arrows follow labeling in the node sequence.}
		\label{Fig_Walks_Node_Traversed_Once}
	\end{center}
\end{figure}

Applying the inclusion exclusion formula in (\ref{Eq_Inclusion_Exclusion_Applied}) produces an exponential number of matrix terms, which do not seem to be solvable for a general hopcount $k$. Instead, we consider here walks in which a node appears exactly once.
\begin{definition}\label{Def_Node_Appearing_Once}
	The set of all possible walks with length $k$, where node $i\in \mathcal{N}$ is traversed exactly once, is denoted as $\mathcal{W}_{(i)}[k]$. The corresponding $N\times N$ walk matrix $M\left(\mathcal{W}_{(i)}[k]\right)$ with the number of such walks between node pair equals
	\[
	\begin{split}
		M\left(\mathcal{W}_{(i)}[k]\right) =& \left(\left(e_{i}\cdot u^{T}\right)\circ A\right)\cdot \left(\left(\left(u-e_{i}\right)\cdot \left(u-e_{i}\right)^{T}\right)\circ A\right)^{k-1} \\
		+& \sum\limits_{m=1}^{k-1} \left(\left(\left(u-e_{i}\right)\cdot \left(u-e_{i}\right)^{T}\right)\circ A\right)^{m-1}\cdot \left(A\circ \left(u\cdot e_{i}^{T}\right)\right) \cdot \left(\left(\left(u-e_{i}\right)\cdot \left(u-e_{i}\right)^{T}\right)\circ A\right)^{k-m} \\
		+& \left(\left(\left(u-e_{i}\right)\cdot \left(u-e_{i}\right)^{T}\right)\circ A\right)^{k-1} \cdot \left(\left(e_{i}\cdot u^{T}\right)\circ A\right).
	\end{split}
	\]
\end{definition}
The $N\times N$ walk matrix $M\left(\mathcal{W}_{(i)}[k]\right)$ consists of $k+1$ terms, because node $i$ can appear on $k+1$ positions in a walk of length $k$. We illustrate examples of walks traversing a node exactly once in Figure \ref{Fig_Walks_Node_Traversed_Once}.

\subsection{Analytic solution for the $N\times N$ path matrix $P_k$}\label{Sec_P_k_Apparing_Once}
Defining the walk sets $\mathcal{W}_{(i)}[k]$ with node $i\in \mathcal{N}$ appearing only once allows us to determine the set of paths $\mathcal{P}[k]$ with hopcount $k$, not by excluding walks with node reappearance from all possible walks $\mathcal{W}[k]$, but instead as the intersection of those walk sets of the form $\mathcal{W}_{(i)}[k], \, i \in \mathcal{N}$
\begin{equation}\label{Eq_Path_Set_Hopcount_k}
	\mathcal{P}[k] = \mathcal{W}_{(i_{0}\in \mathcal{N})}[k] \cap \mathcal{W}_{(i_{1}\in \mathcal{N}\setminus i_{0})}[k] \cap \dots \cap \mathcal{W}_{(i_{k}\in (\mathcal{N}\setminus (i_{0}\cup i_{1}\cup \dots \cup i_{k-1})))}[k].
\end{equation}
\begin{theorem}\label{Theorem_Path_Matrix_Hopcount_k}
	The $N\times N$ path matrix $P_{k}$, whose entries comprise the number of paths with hopcount $k$ between any pair of nodes is defined as follows
	\begin{equation}\label{Eq_Path_Matrix_k}
		P_{k} = \sum_{i_{0}\in \mathcal{N}} \sum_{i_{1}\in \mathcal{N}\setminus i_{0}} \dots \sum_{i_{k}\in \mathcal{N}\setminus (i_{0}\cup i_{1}\cup \dots \cup i_{k-1})} \prod\limits_{z=1}^{k} \left(\left(e_{i_{z-1}}\cdot e_{i_{z}}^{T}\right)\circ A\right),
	\end{equation} or alternatively
	\[
	P_{k} = \sum_{i_{0}= 1}^{N} \sum_{i_{1}= 1}^{N} \dots \sum_{i_{k}= 1}^{N} \prod\limits_{z=1}^{k} \left(\left(e_{i_{z-1}}\cdot e_{i_{z}}^{T}\right)\circ A\right).
	\]
\end{theorem}
\textit{Proof} Relation (\ref{Eq_Path_Matrix_k}) examines all possible labeled sequences of $k+1$ nodes. For each sequence $(i_0,i_1,\dots , i_k)$, we remove all elements from the $N\times N$ adjacency matrix $A$, except for the element $a_{(i_z-1,i_z)}$ between adjacent nodes in the sequence $\left(e_{i_{z-1}}\cdot e_{i_{z}}^{T}\right)\circ A$, where $1\leq z \leq k$. If a node sequence composes a path of lenth $k$, the $(i_0,i_k)$th element of the product $\prod_{z=1}^{k} \left(\left(e_{i_{z-1}}\cdot e_{i_{z}}^{T}\right)\circ A\right)$ equals 1, otherwise 0.
$\hfill \square$.

The solution for the $N\times N$ path matrix $P_{k}$ in (\ref{Eq_Path_Matrix_k}) represents a deterministic counterpart to the relation in \cite[eq. 6]{van2001paths}, defining the probability of a path existence between two nodes. 
For each possible labelled sequence of $k+1$ different nodes, in total $(k+1)!$ of them, relation (\ref{Eq_Path_Matrix_k}) forces all entries of the $N\times N$ adjacency matrix $A$ to zero, expect for the entries between adjacent nodes in the sequence and provides an element one on the position $(ij)$, if the remaining elements compose a path of length $k$ between node $i$ and node $j$. By summing over each possible labelled node sequence, we obtain the $N\times N$ path matrix $P_{k}$, with complexity $O\left( k! k N^3\right)$.

\section{Walks not traversing a node}\label{Sec_Walks_Node_Not_Appearing}

\begin{figure}[!h]
	\begin{center}
		\includegraphics[ angle =0, scale= 0.62]{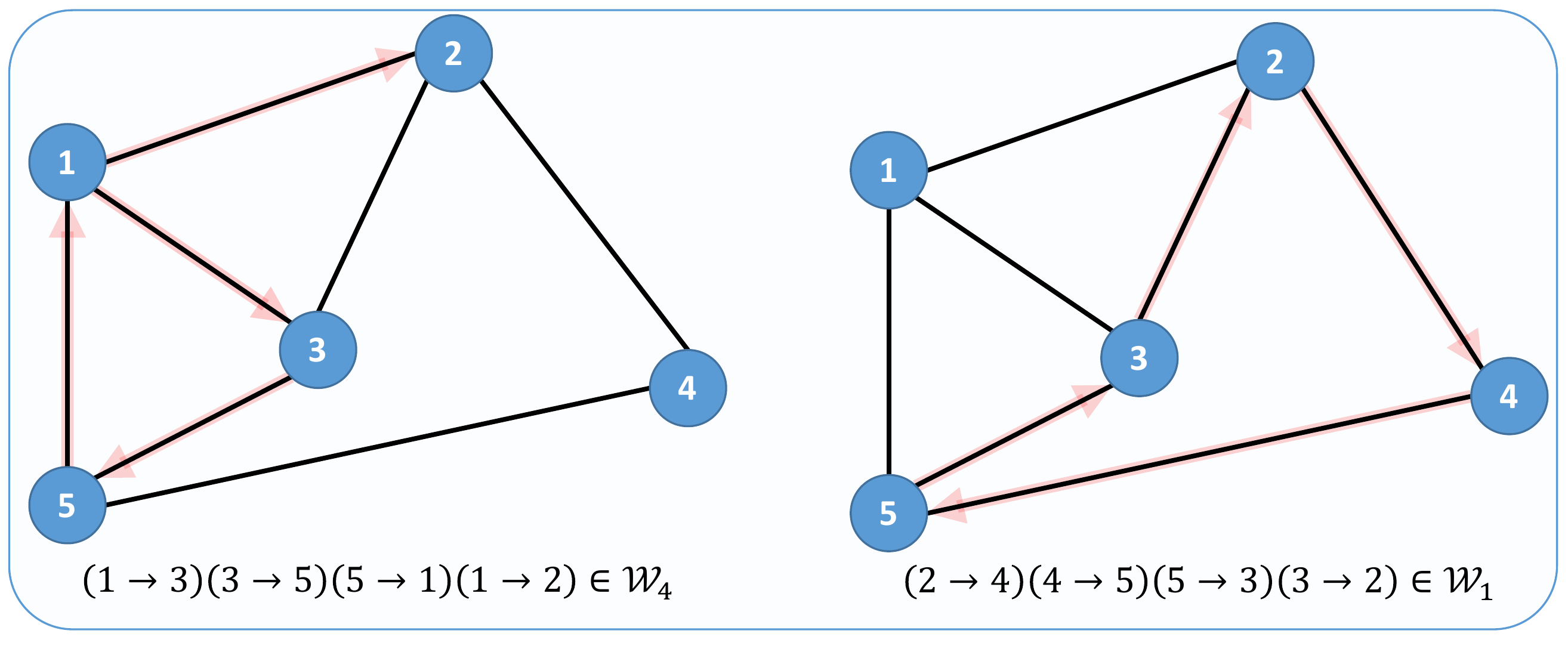}
		\caption{Examples of walks with length $k=N-1$ in which a node is not traversed. Traversed links are colored in red, while arrows follow labeling in the node sequence.}
		\label{Fig_Walks_Node_Not_Traversed}
	\end{center}
\end{figure}

For a general hopcount $k$, there are in total $k!$ matrix terms in (\ref{Eq_Path_Matrix_k}), as each node sequence is labelled. Therefore, an explicit enumeration for an arbitrary hopcount $k$ is infeasible.
However, when computing the $N\times N$ path matrix $P_k$ in a matrix form, labelling nodes in the sequence is not necessary, because the matrix product naturally preserves the information about the source and destination node of each walk, as illustrated in (\ref{Eq_Spliting_Walks}). In this section we introduce walks where a node is not traversed. Originally, this type of walks was defined by Bax in \cite{bax1993inclusion}.
\begin{definition}\label{Def_Walk_Set_Without_Node_m}
	The set of all possible walks with length $k=N-1$, where node $m\in \mathcal{N}$ is not traversed, is denoted as $\mathcal{W}_{m}$. The $N\times N$ corresponding walk matrix with the number of such walks between any pair of nodes equals
	\[
	M\left(\mathcal{W}_{m}\right) = \left(\left(\left(u-e_{m}\right)\cdot \left(u-e_{m}\right)^{T}\right) \circ A\right)^{N-1},
	\]
	where $e_{i}$ denotes the $N\times 1$ basic vector with only one non-zero element $(e_{i})_{i}=1$.
\end{definition}
Figure \ref{Fig_Walks_Node_Not_Traversed} provides two examples of walks not traversing a node.
\subsection{Hamiltonian path matrix $P_{N-1}$}\label{Sec_Hamiltonian_Path_Matrix}
A path of length $N-1$, also known as a Hamiltonian path, is defined by a sequence of $N$ nodes
\[
\left(n_{0},n_{1}, \dots , n_{N-1}\right)
\]
such that $n_{k} \neq n_{l}$, for $0\leq k \leq N-1$ and for $l\in \mathcal{N}\setminus k$, where also $a_{(n_{k},n_{k-1})} = 1$, for $0\leq k \leq N-2$. Because a path by definition consists of different nodes in the sequence, a Hamiltonian path traverses each node in the graph exactly once. Such observation allowed Bax in \cite{bax1993inclusion} to define a set of all walks with hopcount $k=N-1$, where node $m\in \mathcal{N}$ is not traversed. 

Walk sets of the form $\mathcal{W}_{m}$ allow us to define the set of all possible walks of length $k=N-1$ as follows
\begin{equation}\label{Eq_All_Walks_Union_of_sets_without}
	\mathcal{W}[N-1] = \bigcup_{i=1}^{N}\mathcal{W}_{i} \cup \mathcal{P}[N-1]
\end{equation}
Relation (\ref{Eq_All_Walks_Union_of_sets_without}) informs us that either a node is not traversed in a walk of length $N-1$ or that walk represents a path, leading to the following general solution for the $N\times N$ Hamiltonian path matrix $P_{N-1}$
\begin{equation}\label{Eq_Hamiltonian_Path}
	P_{N-1} = A^{N-1} - M\left(\bigcup_{i=1}^{N}\mathcal{W}_{i}\right).
\end{equation}
By applying the inclusion-exclusion formula on the set union from the equation above, Bax obtained
\begin{equation}\label{Eq_Inc_Excl_Solution}
	\begin{split}
		M\left(\bigcup_{i=1}^{N}\mathcal{W}_{i}\right) =& \sum\limits_{i_{1}=1}^{N} M\left(\mathcal{W}_{i_{1}}\right) \\
		-& \sum\limits_{i_{1}=1}^{N-1}\sum\limits_{i_{2}=i_{1}+1}^{N} M\left(\mathcal{W}_{i_{1}}\cap \mathcal{W}_{i_{2}}\right) \\
		+& \sum\limits_{i_{1}=1}^{N-2}\sum\limits_{i_{2}=i_{1}+1}^{N-1}\sum\limits_{i_{3}=i_{2}+1}^{N} M\left(\mathcal{W}_{i_{1}}\cap \mathcal{W}_{i_{2}} \cap \mathcal{W}_{i_{3}}\right) \\
		-& \dots \\
		+& (-1)^{N-3} \sum\limits_{i_{1}=1}^{2}\sum\limits_{i_{2}=i_{1}+1}^{3} \dots \sum\limits_{i_{N}=i_{N-1}+1}^{N} M\left(\bigcap_{z=1}^{N}\mathcal{W}_{i_{z}}\right).
	\end{split}
\end{equation}
A set of intersections from (\ref{Eq_Inc_Excl_Solution}) defines all possible walks with hopcount $N-1$, where multiple nodes are not traversed. The corresponding $N\times N$ walk matrix of such a walk set is
\begin{equation}\label{Eq_Set_Intersections_Without}
	M\left(\mathcal{W}_{i_{1}} \cap \mathcal{W}_{i_{2}} \cap \dots \cap \mathcal{W}_{i_{m}}\right) = \left(\left(\text{diag}\left(u-e_{i_{1}}-e_{i_{2}}-\dots e_{i_{m}}\right)\right)\cdot A \cdot \left(\text{diag}\left(u-e_{i_{1}}-e_{i_{2}}-\dots e_{i_{m}}\right)\right)\right)^{N-1}.
\end{equation}
By combining (\ref{Eq_Hamiltonian_Path}) and (\ref{Eq_Set_Intersections_Without}) Bax derived in \cite{bax1993inclusion} the $N\times N$ Hamiltonian path matrix $P_{N-1}$ as follows
\begin{equation}\label{Eq_Hamiltonian_Paths_Final}
			\small
	\begin{split}
		P_{N-1} =& A^{N-1} - \sum\limits_{i_{1}=1}^{N} \left(\text{diag}\left(u-e_{i_{1}}\right)\cdot A \cdot \text{diag}\left(u-e_{i_{1}}\right) \right)^{N-1} \\
		+& \sum\limits_{i_{1}=1}^{N-1} \sum\limits_{i_{2}=i_{1}+1}^{N} \left(\text{diag}\left(u-e_{i_{1}} - e_{i_{2}}\right)\cdot A \cdot \text{diag}\left(u-e_{i_{1}}-e_{i_{2}}\right) \right)^{N-1} \\
		-& \dots \\
		+& (-1)^{N-3}\cdot \sum\limits_{i_{1}=1}^{2} \sum\limits_{i_{2}=i_{1}+1}^{3} \dots \sum\limits_{i_{N-1}=i_{N-2}+1}^{N} \left(\text{diag}\left(u-\sum\limits_{z=1}^{N-1}e_{i_{z}}\right)\cdot A \cdot \text{diag}\left(u-\sum\limits_{z=1}^{N-1}e_{i_{z}}\right) \right)^{N-1}.
	\end{split}
\end{equation}

\subsection{Analytic solution for the $N\times N$ path matrix $P_k$}\label{Sec_P_k_Not_Apparing}

We here extend the approach of Bax in \cite{bax1993inclusion} and derive an analytic solution for the $N\times N$ path matrix $P_k$ of any hopcount $1\leq k \leq N-1$. The idea behind computing the number of paths with hopcount $k$ between node pairs is to examine all possible unlabeled sequences of $k+1$ nodes, in total $\binom{N}{k+1}$ of them. For each node sequence, we remove links from the graph, not adjacent to any node in the sequence. 
A path of length $k$ in such a reduced graph is equivalent to a Hamiltonian path in the original graph and thus, the idea of Bax from \cite{bax1993inclusion} can be applied.

\begin{theorem}\label{Theorem_Simplified_General_Solution}
    The $N\times N$ path matrix $P_{k}$, whose entries comprise the number of paths with hopcount $k$ between any pair of nodes can be computed as follows
	\begin{equation}\label{Eq_Path_Matrix_k_Simplified}
	\begin{split}
	    &P_k = \sum_{i_0=0}^{N-k-1} \sum_{i_1=i_0 + 1}^{N-k} \dots \sum_{i_k=i_{k-1}+1}^{N} \Bigg[\left( \left(\left(\sum_{z=0}^{k}e_{i_z}\right)\cdot \left(\sum_{z=0}^{k}e_{i_z}^T\right) \right)\circ A\right)^{k} \\
	    -& \sum_{j_0=0}^{k} \left(\left(\left(\sum_{z=0}^{k}e_{i_z}-e_{i_{j_0}}\right) \cdot \left(\sum_{z=0}^{k}e_{i_z}-e_{i_{j_0}}\right)^T\right)\circ A \right)^k \\
	    +& \sum_{j_0=0}^{k - 1} \sum_{j_1=j_0+1}^{k} \left(\left(\left(\sum_{z=0}^{k}e_{i_z}-e_{i_{j_0}}-e_{i_{j_1}}\right) \cdot \left(\sum_{z=0}^{k}e_{i_z}-e_{i_{j_0}}-e_{i_{j_1}}\right)^T\right)\circ A \right)^k \\
	    -& \dots \\
	    +&(-1)^{k-2}\sum_{j_0=0}^{1} \sum_{j_1=j_0+1}^{2} \dots \sum_{j_{k-1}=j_{k-2}+1}^{k} \left(\left(\left(\sum_{z=0}^{k}e_{i_z}-\sum_{q=0}^{k-1}e_{i_{j_q}}\right) \cdot \left(\sum_{z=0}^{k}e_{i_z}-\sum_{q=0}^{k-1}e_{i_{j_q}}\right)^T\right)\circ A \right)^k\Bigg].
	\end{split}
	\end{equation}
\end{theorem}
\textit{Proof} We examine all possible unlabeled sequences of $k+1$ nodes. For each such a sequence $(i_0,i_1,\dots , i_k)$ we transform the $N\times N$ adjacency matrix $A$ by removing all links not adjacent to any node in the sequence $\left(\left(\sum_{z=0}^{k}e_{i_z}\right)\cdot \left(\sum_{z=0}^{k}e_{i_z}^T\right) \right)\circ A$. The modified adjacency matrix allows for applying the inclusion exclusion formula Bax derived in \cite{bax1993inclusion}, because a path of length $k$ in the modified adjacency matrix is equivalent to a Hamiltonian path in the $N\times N$ original adjaceny matrix $A$. $\hfill \square$

There are $\binom{N}{k+1}=\frac{N!}{(k+1)!(N-k-1)!}$ ways to choose $k+1$ nodes out of $N$ nodes. For each set of $k+1$ nodes, relation (\ref{Eq_Path_Matrix_k_Simplified}) defines in total $2^k$ matrix terms\footnote{Each matrix term in (\ref{Eq_Path_Matrix_k_Simplified}) represent the $k$-th power of the adjacency matrix with reduced number of links. Therefore, complexity of computing a matrix term is $O(kN^3)$.} and thus computing the $N\times N$ path matrix $P_k$ implies complexity $O\left(\binom{N}{k+1} k N^3 2^k\right)$.

\subsection{Complexity of computing the $N\times N$ path matrix $P_k$}\label{Sec_Complexity}

\begin{figure}[!h]
	\begin{center}
		\includegraphics[ angle =0, scale= 0.8]{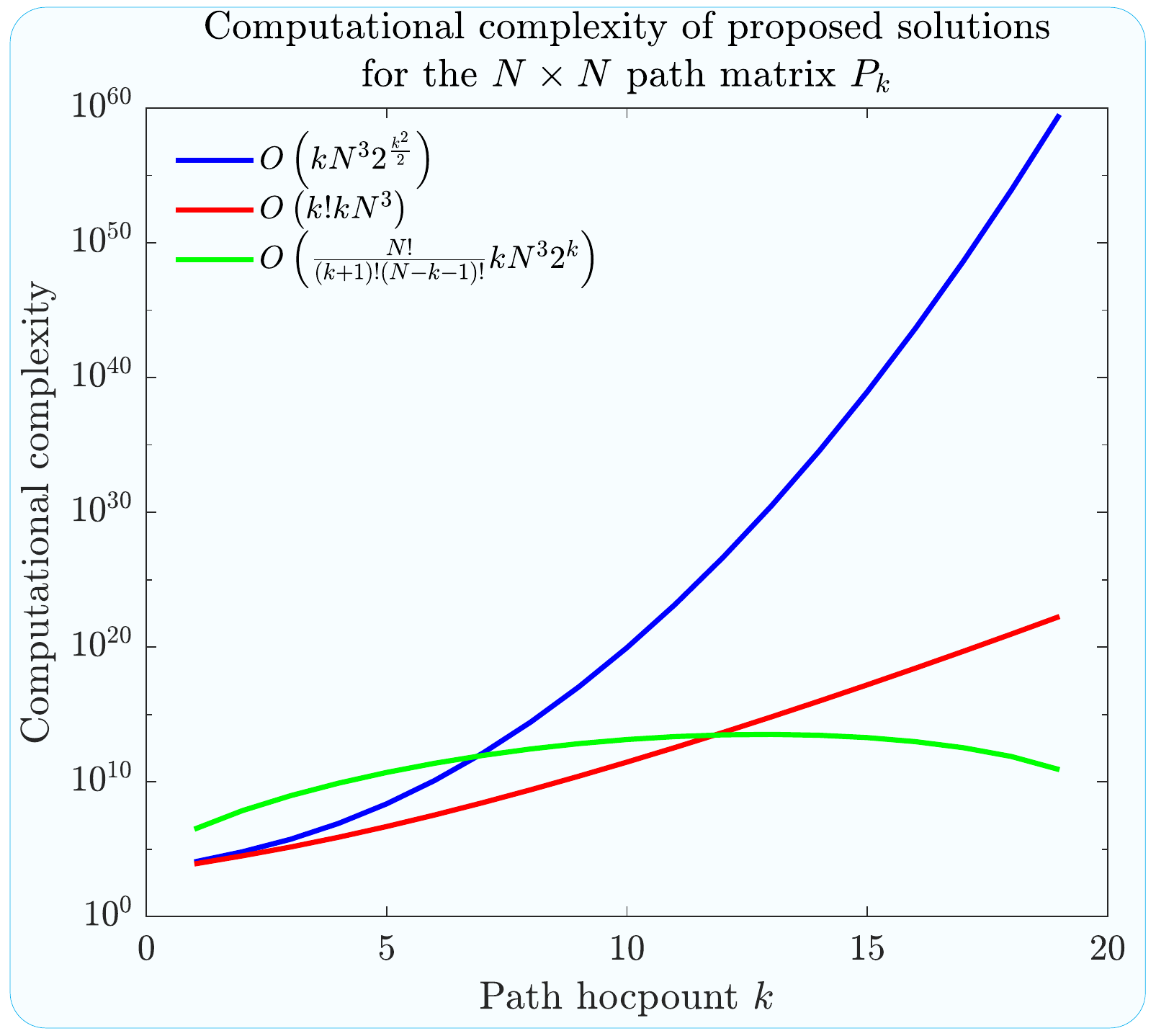}
		\caption{Computational complexity of computing the $N\times N$ path matrix $P_k$ with different hopcount $1\leq k \leq N-1$ using (\ref{Eq_Path_Matrix_General_Sets}) (blue color), (\ref{Eq_Path_Matrix_k}) (red color) and (\ref{Eq_Path_Matrix_k_Simplified})(green color), where $N=20$.}
		\label{Fig_Complexity}
	\end{center}
\end{figure}

We present three analytic solutions for the $N\times N$ path matrix $P_k$, comprising in its entries the number of length $k$ paths between node pairs. 
Figure \ref{Fig_Complexity} provides complexity of computing the $N\times N$ path matrix $P_k$, as a function of the hopcount $k$, for a graph of $N=20$ nodes. 
Complexity of the solution in (\ref{Eq_Path_Matrix_General_Sets}) (blue color), based on walks traversing a node multiple times, represents the most complex approach for almost the entire range of hopcount $k$ values. Despite its complexity, for small $k$, the solution in (\ref{Eq_Path_Matrix_General_Sets}) is the most insightful, from a linear algebra point of view.

Computing the $N\times N$ path matrix $P_k$ using (\ref{Eq_Path_Matrix_k}) (red color in Figure \ref{Fig_Complexity}) requires the least computational effort, for smaller values of the hopcount $k$, because it examines all possible labeled sequences of $k+1$ nodes. In contrast, for smaller values of hopcount $k$, the third approach in  (\ref{Eq_Path_Matrix_k_Simplified}) (presented in green color in Figure \ref{Fig_Complexity})) is far more computationally demanding, because it applies the inclusion exclusion formula on each possible unlabeled sequence of $k+1$ nodes. As the path length $k$ increases, there are less unlabeled sequences of nodes, allowing the third solution in (\ref{Eq_Path_Matrix_k_Simplified}) to perform the best, in terms of complexity.

\section{Recursive algorithm for computing the number of paths}\label{Sec_Recursive_Algorithm}

We provide a simple recursive algorithm \ref{Fig_Alg_Recursion} that computes the $N\times N$ path matrices $P_{k}$, where $1\leq k \leq N-1$. The proposed algorithm treats each path independently (via recursions) and prevents repeating nodes in the walk sequence.

\begin{figure}[ptbh]
	\begin{center}
		\framebox[4.5in]{\begin{minipage}{4in}
				{\sc DeterminePaths} ($A, \, N$)
				\begin{tabbing}
					\bf{Input}: $A, \, N$\\
					\bf{Output}: $P_{1}, \, P_{2}, \, \dots , \, P_{N-1}$\\
					1. Initialise node-based path matrix $T \gets O_{N-1\times N}$\\ 
					2. Initialise path matrices $P_{k} \gets O_{N\times N}$, where $1\leq k < N$\\ 
					3. \bf{for} \= $i \gets 1$ to $N$ \\
					4. \> $T \gets $ {\sc ComputePaths} ($O_{N\times N}, \, A, \, i, \, i, \, 0$) \\
					5. \> \bf{for} \= $j \gets 1$ to $N-1$ \\
					6. \> \> Store $j$-th row of $T$ as the $i$-th row of $P_{j}$\\
					7. \> \bf{end for} \\
					8. \bf{end for} \\
					9. \bf{return} $P_{1}, \, P_{2}, \, \dots , \, P_{N-1}$ \\
				\end{tabbing}
		\end{minipage}}
	\end{center}
	\caption{Pseudocode for calling the recursive Algorithm for determining all paths in a graph, with the network size $N$ and the the $N\times N$ adjacency matrix $A$ as input.}%
	\label{Fig_Alg_Calling}%
\end{figure}

The proposed recursive algorithm identifies each possible path in a graph and increment the corresponding element of the $N\times N$ path matrix $P_{k}$, with $1\leq k<N$. For each node $i\in \mathcal{N}$ in $G$, we set the hopcount to $0$ and call the recursive procedure, as provided in line $4$ of the pseudocode \ref{Fig_Alg_Calling}. The recursive algorithm \ref{Fig_Alg_Recursion} returns the $N-1\times N$ node $i$ based path matrix $T$, where the element $T_{jm}$ denotes the number of length $j$ paths between node $i$ and node $m$. Therefore, in line $6$ we store the $j$-th row of the $N-1\times N$ node based path matrix $T$ as the $i$-th row of the $N\times N$ path matrix $P_j$. 

\begin{figure}[ptbh]
	\begin{center}
		\framebox[6.5in]{\begin{minipage}{4in}
				{\sc ComputePaths} ($T, \, A, \, n_{0}, \, n_{k}, \, k$)
				\begin{tabbing}
					\bf{Input}: $T, \, A, \, n_{0}, \, n_{k}, \, k$\\
					\bf{Output}: $T, \, A, \, n_{0}, \, n_{k}, \, k$\\
					1. Increment the hopcount $k \gets k+1$ \\
					2. Determine the destination node $n_{k}$ neighbours $\mathcal{N}_{n_{k}} \gets \{j\,| a_{n_{k},j} = 1, j \in \mathcal{N}\}$\\
					3. Increment the node-based path matrix $T_{k,j} \gets T_{k,j} + 1$, where $j \in \mathcal{N}_{n_{k}}$ \\
					4. Delete links adjacent to the destination node $a_{n_{k},j} \gets 0$ and $a_{j,n_{k}} \gets 0$, where $j \in \mathcal{N}_{n_{k}}$  \\
					5. \bf{if} \= $|\mathcal{N}_{n_{k}}| = 0$ or $k = N - 1$ \\
					6. \> \bf{return} \= ($T, \, A, \, n_{0}, \, n_{k}, \, k$) \\
					7. \bf{else} \\
					8. \> $\hat{k} \gets k$ and $\hat{A} \gets A$ \\
					9. \> \bf{for} \= $m \gets 1$ to $|\mathcal{N}_{n_{k}}|$ \\
					10. \> \> \bf{return} ($T, \, A, \, n_{0}, \, j_{m}, \, k$) $\gets$ {\sc ComputePaths} ($T, \, \hat{A}, \, n_{0}, \, j_{m}, \, \hat{k}$) \\
					11. \> \bf{end for} \\
					12. \bf{end if} \\
				\end{tabbing}
		\end{minipage}}
	\end{center}
	\caption{Metacode of the Algorithm for determining all paths in a graph, with
		the $N\times N$ adjacency matrix $A$, the $N\times N$ node-based path matrix $T$, source node $n_{0}$, destination node $n_{k}$ and hopcount $k$ as input.}%
	\label{Fig_Alg_Recursion}%
\end{figure}

The recursive procedure in Algorithm \ref{Fig_Alg_Recursion} receives as input the source node $n_{0}$ and the destination node $n_{k}$, as well the $N\times N$ adjacency matrix $A$ and the current hopcount $k$. We firstly identify neighbours of the destination node $j \in \mathcal{N}_{n_{k}}$ in line $2$. After incrementing the hopcount $k$ in line $3$, in the following step we account for the paths reaching any of the neighbours in $\mathcal{N}_{n_{k}}$. A crucial step is to remove all links adjacent to the destination node $n_{k}$, defined in line $4$, as not to allow paths to reach again node $n_{k}$. Upon removing links, for each neighbour $j$, we call the recursive algorithm in line $10$, with the incremented hopcount and the updated destination node $j$. Therefore, within the recursive function, all links that would lead to node reappearance are removed, allowing us to determine each adjacent node of the destination node as a valid path extension. The recursion stops once a destination node has no neighbours or the hopcount equals $k=N-1$, as defined in line $5$.

\begin{figure}[!h]
	\begin{center}
		\includegraphics[ angle =0, scale= 0.58]{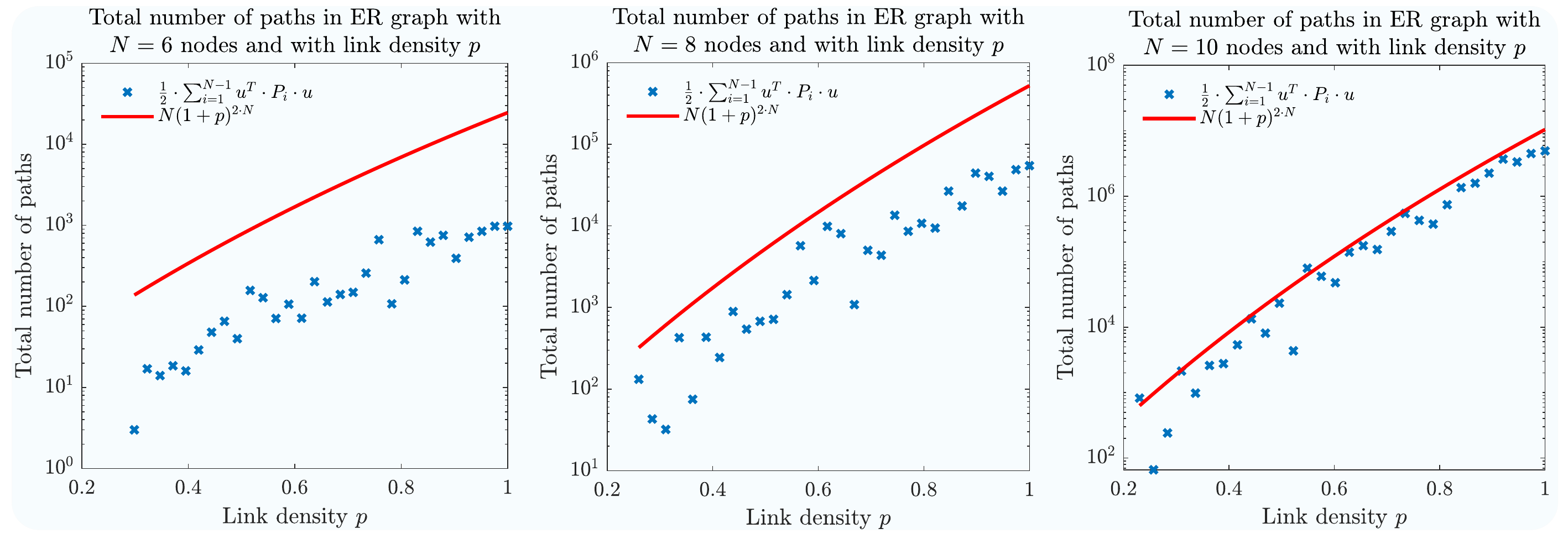}
		\caption{Total number of paths in Erdős Rényi graph with $N=6$ (left figure), $N=8$ nodes (middle figure) and $N=10$ nodes (right figure), for different values of the link density $p$.}
		\label{Fig_ER_N_10_All_Paths_per_p}
	\end{center}
\end{figure}

By performing the Algorithm \ref{Fig_Alg_Recursion} once, we obtain information about each possible path in a graph. Since the proposed recursive algorithm \ref{Fig_Alg_Recursion} accounts for each possible path, its complexity scales linearly with the total number of paths $\frac{1}{2}\cdot \sum_{i=1}^{N-1}u^{T}\cdot P_{i}\cdot u$. Therefore, Figure \ref{Fig_ER_N_10_All_Paths_per_p} depicts the total number of paths of Erdős Rényi (ER) random graph with $N=10$ nodes, for different link density $p$, revealing an exponential correlation. The complexity of the proposed algorithm $O(N(1+p)^{2N})$ scales exponentially with the network size $N$. Number of paths with hopcount $k$ between node pairs can be determined with the proposed recursive algorithm with complexity $O(N(1+p)^{2k})$.

\section{Conclusion}\label{Sec_Conclusion}
We introduce three types of walks: walks with a node reappearing in the sequence, walks traversing a node exactly once, and those not traversing a node. Based on considered walk types, we derive analytic solutions for the number of paths of a certain length between node pairs in a matrix form. Depending on the path length, different solutions require the least computational effort.
We propose a recursive algorithm for determining all possible paths between node pairs, whose complexity scales linearly with the total number of paths in a graph. The proposed recursive algorithm applies to a directed (un)weighted network.

\bibliography{20220919_Network_Paths_ArXiv.bib}{}
\bibliographystyle{ieeetr}

\appendix
 
\section{$N\times N$ matrix $F_{k}$}\label{App_Theorem_Recursive}

We derive the first two sum terms of the $N\times N$ matrix $F_k$ in (\ref{Eq_Recursive_Solution_Main}). Firstly, we split the first term in (\ref{Eq_Inclusion_Exclusion_Applied}) into two sum terms, where for the second sum term it holds $j_{2}=k$
\[
\sum\limits_{i_1=0}^{k-2}\sum\limits_{j_1=i+2}^{k} M\left(\mathcal{W}_{(i_1,j_1)}[k]\right) = \sum\limits_{i_{1}=0}^{k-3}\sum\limits_{j_{1}=i_{1}+2}^{k-1} M\left(\mathcal{W}_{(i_{1},j_{1})}[k]\right) + \sum\limits_{i_{2}=0}^{k-2} M\left(\mathcal{W}_{(i_{2},k)}[k]\right).
\]
Because the counters $i_1<k$ and $j_1<k$ of the first sum terms in the equation above are always smaller than $k$, we import (\ref{Eq_Repeating_Node_Walk_Matrix}) and further transform the equation above
\[
\sum\limits_{i_1=0}^{k-2}\sum\limits_{j_1=i+2}^{k} M\left(\mathcal{W}_{(i_1,j_1)}[k]\right) - \sum\limits_{i_{1}=0}^{k-3}\sum\limits_{j_{1}=i_{1}+2}^{k-1} M\left(\mathcal{W}_{(i_{1},j_{1})}[k-1]\right)\cdot A = \sum\limits_{i_{2}=0}^{k-2} A^{i_2}\cdot \left(I\circ A^{k-i_2}\right).
\]
Similarly, we transform the second sum term in (\ref{Eq_Inclusion_Exclusion_Applied}) by excluding the case $j_2 = k$ from the sum term
\begin{equation}\label{Eq_Proof_Recursive_Term_2}
\resizebox{.91 \textwidth}{!}{$
	\begin{split}
		\sum\limits_{i_{1}=0}^{k-2}\sum\limits_{j_{1}=i_{1}+2}^{k}\sum\limits_{i_{2}=i_{1}}^{k-2}\sum\limits_{j_{2}=q_2}^{k} M\left(\mathcal{W}_{(i_{1},j_{1})}[k] \cap \mathcal{W}_{(i_{2},j_{2})}[k]\right) & = 	\sum\limits_{i_{1}=0}^{k-2}\sum\limits_{j_{1}=i_{1}+2}^{k}\sum\limits_{i_{2}=i_{1}}^{k-2}\sum\limits_{j_{2}=q_2}^{k-1} M\left(\mathcal{W}_{(i_{1},j_{1})}[k] \cap \mathcal{W}_{(i_{2},j_{2})}[k]\right) \\
		& + \sum\limits_{i_{1}=0}^{k-2}\sum\limits_{j_{1}=i_{1}+2}^{k}\sum\limits_{i_{2}=i_{1}}^{k-2} M\left(\mathcal{W}_{(i_{1},j_{1})}[k] \cap \mathcal{W}_{(i_{2},k)}[k]\right)
	\end{split}$}
\end{equation}

We further exclude the case $j_{1}=k$ from both sum term on the right-hand side of (\ref{Eq_Proof_Recursive_Term_2}) and obtain

\begin{equation}\label{Eq_Proof_Recursive_Term_Full}
\resizebox{.91 \textwidth}{!}{$
	\begin{split}
		\sum\limits_{i_{1}=0}^{k-2}\sum\limits_{j_{1}=i_{1}+2}^{k}\sum\limits_{i_{2}=i_{1}}^{k-2}\sum\limits_{j_{2}=q_2}^{k} M\left(\mathcal{W}_{(i_{1},j_{1})}[k] \cap \mathcal{W}_{(i_{2},j_{2})}[k]\right) & = 	\sum\limits_{i_{1}=0}^{k-2}\sum\limits_{j_{1}=i_{1}+2}^{k-1}\sum\limits_{i_{2}=i_{1}}^{k-2}\sum\limits_{j_{2}=q_2}^{k-1} M\left(\mathcal{W}_{(i_{1},j_{1})}[k] \cap \mathcal{W}_{(i_{2},j_{2})}[k]\right) \\
		& + \sum\limits_{i_{1}=0}^{k-2}\sum\limits_{i_{2}=i_{1}}^{k-2}\sum\limits_{j_{2}=q_2}^{k-1} M\left(\mathcal{W}_{(i_{1},k)}[k] \cap \mathcal{W}_{(i_{2},j_{2})}[k]\right)\\
		& + \sum\limits_{i_{1}=0}^{k-2}\sum\limits_{j_{1}=i_{1}+2}^{k-1}\sum\limits_{i_{2}=i_{1}}^{k-2} M\left(\mathcal{W}_{(i_{1},j_{1})}[k] \cap \mathcal{W}_{(i_{2},k)}[k]\right) \\ 
		& + \sum\limits_{i_{1}=0}^{k-2}\sum\limits_{i_{2}=i_{1}}^{k-2} M\left(\mathcal{W}_{(i_{1},k)}[k] \cap \mathcal{W}_{(i_{2},k)}[k]\right)
	\end{split}$}
\end{equation}

The first sum term can be split in two sub walks as in (\ref{Eq_Spliting_Walks}), transforming (\ref{Eq_Proof_Recursive_Term_Full}) as follows

\begin{equation}\label{Eq_Proof_Recursive_Term_Optimised}
\resizebox{.91 \textwidth}{!}{$
	\begin{split}
		\sum\limits_{i_{1}=0}^{k-2}\sum\limits_{j_{1}=i_{1}+2}^{k}\sum\limits_{i_{2}=i_{1}}^{k-2}\sum\limits_{j_{2}=q_2}^{k} M\left(\mathcal{W}_{(i_{1},j_{1})}[k] \cap \mathcal{W}_{(i_{2},j_{2})}[k]\right) & = 	\sum\limits_{i_{1}=0}^{k-3}\sum\limits_{j_{1}=i_{1}+2}^{k-1}\sum\limits_{i_{2}=i_{1}}^{k-3}\sum\limits_{j_{2}=q_2}^{k-1} M\left(\mathcal{W}_{(i_{1},j_{1})}[k-1] \cap \mathcal{W}_{(i_{2},j_{2})}[k-1]\right)\cdot A \\
		& +  \sum\limits_{i_{1}=0}^{k-2}\sum\limits_{i_{2}=i_{1}}^{k-2}\sum\limits_{j_{2}=q_2}^{k-1} A^{i_1}\cdot \left(I\circ \left(A^{i_{2}-i_{1}}\cdot \left(I\circ A^{j_{2}-i_{2}}\right)\cdot A^{k-j_{2}}\right)\right)\\
		& + \sum\limits_{i_{1}=0}^{k-2}\sum\limits_{j_{1}=i_{1}+2}^{k-1}\sum\limits_{i_{2}=i_{1}}^{k-2} M\left(\mathcal{W}_{(i_{1},j_{1})}[k] \cap \mathcal{W}_{(i_{2},k)}[k]\right) \\ 
		& + \sum\limits_{i_{1}=0}^{k-4}\sum\limits_{i_{2}=i_{1}+2}^{k-2} A^{i_{1}}\cdot \left(I\circ A^{i_{2}-i_{1}}\right)\cdot \left(I\circ A^{k-i_{2}}\right)
	\end{split}$}
\end{equation}
The third sum term on the right-hand side of the relation above can be split into three sum terms, where $i_{2}<j_{1}$, secondly $i_{2}=j_{1}$ and finally where $i_{2} > j_{1}$, allowing us to transform (\ref{Eq_Proof_Recursive_Term_Optimised}) as follows
\begin{equation}\label{Eq_Proof_Recursive_Term_Extended}
\resizebox{.91 \textwidth}{!}{$
	\begin{split}
		\sum\limits_{i_{1}=0}^{k-2}\sum\limits_{j_{1}=i_{1}+2}^{k}\sum\limits_{i_{2}=i_{1}}^{k-2}\sum\limits_{j_{2}=q_2}^{k} M\left(\mathcal{W}_{(i_{1},j_{1})}[k] \cap \mathcal{W}_{(i_{2},j_{2})}[k]\right) & = 	\sum\limits_{i_{1}=0}^{k-3}\sum\limits_{j_{1}=i_{1}+2}^{k-1}\sum\limits_{i_{2}=i_{1}}^{k-3}\sum\limits_{j_{2}=q_2}^{k-1} M\left(\mathcal{W}_{(i_{1},j_{1})}[k-1] \cap \mathcal{W}_{(i_{2},j_{2})}[k-1]\right)\cdot A \\
		& +  \sum\limits_{i_{1}=0}^{k-2}\sum\limits_{i_{2}=i_{1}}^{k-2}\sum\limits_{j_{2}=q_2}^{k-1} A^{i_1}\cdot \left(I\circ \left(A^{i_{2}-i_{1}}\cdot \left(I\circ A^{j_{2}-i_{2}}\right)\cdot A^{k-j_{2}}\right)\right)\\
		& + \sum\limits_{i_{1}=0}^{k-2}\sum\limits_{j_{1}=i_{1}+2}^{k-1}\sum\limits_{i_{2}=i_{1}}^{j_{1}-1} M\left(\mathcal{W}_{(i_{1},j_{1})}[k] \cap \mathcal{W}_{(i_{2},k)}[k]\right) \\ 
		& + \sum\limits_{i_{1}=0}^{k-2}\sum\limits_{j_{1}=i_{1}+2}^{k-1} M\left(\mathcal{W}_{(i_{1},j_{1})}[k] \cap \mathcal{W}_{(j_{1},k)}[k]\right) \\ 
		& + \sum\limits_{i_{1}=0}^{k-2}\sum\limits_{j_{1}=j_{1}+2}^{k-1}\sum\limits_{i_{2}=i_{1}}^{k-2} M\left(\mathcal{W}_{(i_{1},j_{1})}[k] \cap \mathcal{W}_{(i_{2},k)}[k]\right) \\ 
		& + \sum\limits_{i_{1}=0}^{k-4}\sum\limits_{i_{2}=i_{1}+2}^{k-2} A^{i_{1}}\cdot \left(I\circ A^{i_{2}-i_{1}}\right)\cdot \left(I\circ A^{k-i_{2}}\right).
	\end{split}$}
\end{equation}
Finally, the second sum term of the $N\times N$ matrix $F_k$ is defined as follows
\begin{equation}\label{Eq_Proof_Recursive_Term_2_Final}
\resizebox{.91 \textwidth}{!}{$
	\begin{split}
		\sum\limits_{i_{1}=0}^{k-2}\sum\limits_{j_{1}=i_{1}+2}^{k}\sum\limits_{i_{2}=i_{1}}^{k-2}\sum\limits_{j_{2}=q_2}^{k} M\left(\mathcal{W}_{(i_{1},j_{1})}[k] \cap \mathcal{W}_{(i_{2},j_{2})}[k]\right) & = 	\sum\limits_{i_{1}=0}^{k-3}\sum\limits_{j_{1}=i_{1}+2}^{k-1}\sum\limits_{i_{2}=i_{1}}^{k-3}\sum\limits_{j_{2}=q_2}^{k-1} M\left(\mathcal{W}_{(i_{1},j_{1})}[k-1] \cap \mathcal{W}_{(i_{2},j_{2})}[k-1]\right)\cdot A \\
		& +  \sum\limits_{i_{1}=0}^{k-2}\sum\limits_{i_{2}=i_{1}}^{k-2}\sum\limits_{j_{2}=q_2}^{k-1} A^{i_1}\cdot \left(I\circ \left(A^{i_{2}-i_{1}}\cdot \left(I\circ A^{j_{2}-i_{2}}\right)\cdot A^{k-j_{2}}\right)\right)\\
		& + \sum\limits_{i_{1}=0}^{k-2}\sum\limits_{j_{1}=i_{1}+2}^{k-1}\sum\limits_{i_{2}=i_{1}}^{j_{1}-1} A^{i_{1}}\cdot \left(A^{i_{2}-i_{1}}\circ A^{j_{1}-i_{2}}\circ A^{k-j_{1}}\right)  \\ 
		& + \sum\limits_{i_{1}=0}^{k-4}\sum\limits_{j_{1}=i_{1}+2}^{k-2} A^{i_{1}}\left(I\circ A^{j_{1}-i_{1}}\right)\cdot\left(I\circ A^{k-j_{1}}\right) \\ 
		& + \sum\limits_{i_{1}=0}^{k-2}\sum\limits_{j_{1}=i_{1}+2}^{k-1}\sum\limits_{i_{2}=j_{1}+1}^{k-2} A^{i_{1}}\cdot \left(I\circ A^{j_{1}-i_{1}}\right)\cdot A^{i_{2}-j_{1}}\cdot \left(I\circ A^{k-i_{2}}\right) \\ 
		& + \sum\limits_{i_{1}=0}^{k-4}\sum\limits_{i_{2}=i_{1}+2}^{k-2} A^{i_{1}}\cdot \left(I\circ A^{i_{2}-i_{1}}\right)\cdot \left(I\circ A^{k-i_{2}}\right).
	\end{split}$}
\end{equation}

\end{document}